\documentclass[10pt]{amsart}
\usepackage{mathptmx}
\usepackage{amsmath}     
\usepackage{amssymb}
\usepackage{array}
\usepackage{geometry}
\usepackage[bookmarks=true,colorlinks=true, pdfstartview=FitV, linkcolor=black, citecolor=blue, urlcolor=black]{hyperref}
\usepackage{movie15}

\usepackage{color}
\definecolor{DarkRed}{rgb}{0.55,.00,0.2}
\definecolor{DarkGrey}{rgb}{0.35,.35,0.35}

\theoremstyle{definition}

\theoremstyle{remark}

\numberwithin{equation}{section}

\begin{document}

\title{Orthogonal  polynomials   with the Prudnikov-type weights}

\author{S. Yakubovich}
\address{Department of Mathematics, Fac. Sciences of University of Porto,Rua do Campo Alegre,  687; 4169-007 Porto (Portugal)}
\email{ syakubov@fc.up.pt}
\thanks{\thanks{The work was partially supported by CMUP [UID/MAT/00144/2019], which is funded by FCT(Portugal) with national (MEC) and European structural funds through the programs FEDER, under the partnership agreement PT2020, and Project STRIDE - NORTE-01-0145-FEDER- 000033, funded by ERDF - NORTE 2020.  The author thanks Marco Martins Afonso for  necessary numerical calculations and verifications of some formulas. 
 }}


\date{\today}


\begin{abstract}  New sequences of orthogonal polynomials  with respect to the weight functions $e^{-x} \rho_\nu(x),\  e^{- 1/x}  x^{-1} \rho_{\nu} (x), \\\rho_{\nu}(x)= 2 x^{\nu/2} K_\nu(2\sqrt x),\ x >0, \nu \in \mathbb{R}$, where  $K_\nu(z)$ is the modified Bessel function, are investigated.  The recurrence relations, explicit representations, generating functions and Rodrigues-type formulae  are obtained.   \\
\end{abstract}

\maketitle

\markboth{\rm \centerline{ S. Yakubovich}}{Orthogonal  polynomials  with the Prudnikov-type weights}

\section{Introduction and preliminary results}

In the preceding paper \cite{YAP} the author gave an interpretation of the Prudnikov sequence of orthogonal polynomials  with  the weight function $x^\alpha \rho_\nu(x)$, where \  $ \rho_{\nu}(x)= 2 x^{\nu/2} K_\nu(2\sqrt x),\ x >0, \nu \ge 0,\ \alpha > -1$ and  $K_\nu(z)$ is  the modified Bessel function  \cite{Bateman}, Vol. II, in terms of the so-called composition orthogonality with respect to the measure $x^\nu e^{-x} dx$ related to the classical associated Laguerre polynomials \cite{Chi}.  The object of the present paper is to  characterize sequences of orthogonal polynomials, satisfying the following orthogonality conditions

$$\int_{0}^{\infty}  P^\nu_{m} (x) P^\nu_{n}(x) e^{-x} \rho_\nu(x) dx = \delta_{n,m},\eqno(1.1)$$

$$\int_{0}^{\infty}  Q^\nu_{m} (x) Q^\nu_{n}(x) e^{- 1/x} \rho_\nu(x) {dx\over x} = \delta_{n,m},\eqno(1.2)$$
where $\delta_{n,m},\   n,m\in\mathbb{N}_{0}$ is the Kronecker symbol.   Namely, our goal will be  explicit expressions of these sequences, the Rodrigues-type formulas, the 3-term recurrence relations and generating functions associated with these orthogonal polynomials.  Conditions (1.1), (1.2) are equivalent to the following equalities

$$\int_{0}^{\infty}  P^\nu_{n} (x)  e^{-x} \rho_\nu(x) \  x^m dx = 0,\quad m=0,1,\dots, n-1,\eqno(1.3)$$

$$\int_{0}^{\infty}  Q^\nu_{n} (x)  e^{- 1/x} \rho_\nu(x)\ x^{m-1} dx =  0,\quad   m=0,1,\dots, n-1. \eqno(1.4)$$
The weight  function  $\rho_{\nu}(x)$  has  the Mellin-Barnes integral representation in the form (cf. \cite{YAP})
$$
\rho_\nu(x)=  \frac{1}{2\pi i} \int_{\gamma-i\infty}^{\gamma+i\infty} \Gamma(\nu+s) \Gamma (s) x^{-s} ds\ , \quad x, \gamma \in \mathbb{R}_{+},\ \nu  \in \mathbb{R},\eqno(1.5)
$$
where $\Gamma(z)$ is Euler's gamma-function \cite{Bateman}, Vol. I. Moreover, the Parseval equality for the Mellin transform \cite{Tit} permits to write (1.5) as the Laplace integral representation for this weight function. In fact, we obtain 

$$\rho_\nu(x)= \int_0^\infty t^{\nu-1} e^{-t - x/t} dt,\quad  x >0,\ \nu \in \mathbb{R}.\eqno(1.6)$$ 
The asymptotic behavior of the modified Bessel function at infinity and near the origin \cite{Bateman}, Vol. II gives the corresponding values for  function $\rho_\nu(x),\ \nu \in \mathbb{R}$.  Precisely, we have
$$\rho_\nu (x)= O\left( x^{(\nu-|\nu|)/2}\right),\  x \to 0,\ \nu\neq 0, \quad  \rho_0(x)= O( \log x),\ x \to 0,\eqno(1.7)$$

$$ \rho_\nu(x)= O\left( x^{\nu/2- 1/4} e^{- 2\sqrt x} \right),\ x \to +\infty.\eqno(1.8)$$
The moments of the weights $\omega^+_\nu(x)= e^{-x} \rho_\nu(x),\    \omega^-_\nu(x) = e^{- 1/x}  x^{-1} \rho_{\nu} (x)$ can be calculated as follows. Precisely, taking (1.6) and changing the order of integration by Fubini's theorem, we deduce after calculation the classical Euler integral 

$$\int_0^\infty e^{-x}  \rho_\nu(x) x^\mu dx =  \int_0^\infty e^{-x}   x^\mu  \int_0^\infty t^{\nu-1} e^{-t - x/t} dt dx $$

$$=  \Gamma(1+\mu)  \int_0^\infty {t^{\nu+\mu}  e^{-t } \over (1+t)^{\mu+1} }  dt,\quad  \mu +\nu,\ \mu  > -1. $$
But the latter integral with respect to $t$ can be expressed in terms of the Tricomi function $\Psi(a,b; z)$ \cite{Bateman},  Vol. I and \cite{PrudnikovMarichev}, Vol. I, Entry 2.3.6.9.  Hence we get the formula for the  moments of  $\omega^+_\nu(x)$

$$ \int_0^\infty \omega^+_\nu(x) x^\mu dx = \Gamma(\mu+\nu+1)\Gamma(\mu+1) \Psi \left(1+\mu+\nu,\ 1+\nu; 1\right).\eqno(1.8)$$
Analogously, the moments of $ \omega^-_\nu(x)$ can be obtained, employing relation (2.16.8.13) in \cite{PrudnikovMarichev}, Vol. II.  Thus it gives 

$$ \int_0^\infty \omega^-_\nu(x) x^{\mu} dx =  4 \int_0^\infty e^{-1/x^2}   x^{2\mu+\nu -1} K_\nu(2 x) dx = \Gamma(\mu) \Gamma(\mu+\nu) \ 
{}_0F_2\left(1- \mu-\nu, \ 1-\mu;\ -1\right)  $$

$$+ \Gamma(- \nu) \Gamma(- \mu-\nu) \  {}_0F_2\left(1+ \mu+ \nu, \ 1+\nu;\ -1\right)  + \Gamma( \nu) \Gamma(- \mu) \  {}_0F_2\left(1+ \mu, \ 1-\nu;\ -1\right),\eqno(1.9)$$
where ${}_0F_2(a,b; z)$ is the generalized hypergeometric function \cite{Bateman}, Vol. I.

Further, denoting the operator of the so-called Laguerre derivative by $\beta = DxD$ and its companion $\theta =xDx$  \cite{YAP}, where $D$ is the differential operator $D= {d\over dx}$, we calculate them $n$-th power, appealing to the  Viskov-type identities \cite{Viskov}

$$\beta^n = \left( DxD\right)^n = D^n x^n D^n,\quad  \theta^n = \left( xDx\right)^n = x^n D^n x^n,\quad   n \in\mathbb{N}_{0}.\eqno(1.10)$$
 In particular, we easily find the formulas

$$ (\beta^n  \rho_0) (x)= \left( DxD\right)^n  \rho_0 = \rho_0(x), \quad   (\beta^n  \rho_1) (x) = \left( DxD\right)^n  \rho_1 = \rho_1(x) - n \rho_0(x),\quad   n \in\mathbb{N}_{0},\eqno(1.11)$$
 
 $$ (\theta^n  1) (x) = \left( xDx\right)^n  1= n! x^n,\quad  (\theta^n  x^k) (x)=  \left( xDx\right)^n  x^k = {(n+k)!\over k!} x^{n+k},  \quad  n, k  \in\mathbb{N}_{0}.\eqno(1.12)$$
The quotient of the scaled Macdonald functions $\rho_\nu, \rho_{\nu+1}$ is given by the important Ismail integral representation \cite{Ismail}

$${\rho_\nu(x) \over \rho_{\nu+1}(x) } = {1\over \pi^2} \int_0^\infty {y^{-1} dy \over (x+y) \left[ J_{\nu+1}^2 (2\sqrt y)+ Y^2_{\nu+1}(2\sqrt y) \right] },\eqno(1.13)$$
where $J_\nu(z), Y_\nu(z)$ are Bessel functions of the first and second kind, respectively \cite{Bateman}.  Meanwhile, Entry 2.19.4.13  in \cite{PrudnikovMarichev}, Vol. II represents the product $x^n \rho_\nu(x)$ in terms of the associated Laguerre polynomials \cite{Chi}.  Precisely, it has

$${(-1)^n x^n\over n!}\  \rho_\nu(x)=   \int_0^\infty t^{\nu+n -1} e^{-t - x/t}  L_n^\nu(t) dt,\quad    n \in\mathbb{N}_{0}.\eqno(1.14)$$
Moreover,  the Riemann-Liouville fractional integral \cite{YaL}

$$ \left( I_{-}^\alpha  f \right) (x)  = {1\over \Gamma(\alpha)} \int_x^\infty (t-x)^{\alpha-1} f(t) dt,\ {\rm Re} \alpha > 0,\eqno(1.15)$$
and Entry 2.16.3.8  in \cite{Bateman}, Vol. II drive to the  formula 
$$\rho_\nu(x)= \left( I_{-}^\nu \rho_0 \right) (x).\eqno(1.16)$$
Moreover, the index law for fractional integrals immediately implies

$$ \rho_{\nu+\mu} (x)= \left( I_{-}^\nu \rho_\mu \right) (x)=   \left( I_{-}^\mu \rho_\nu \right) (x).\eqno(1.17)$$
The corresponding definition of the fractional derivative presumes the relation $ D^\mu_{-}= - D  I_{-}^{1-\mu}$.   Hence for the ordinary $n$-th derivative of $\rho_\nu(x)$ we get

$$D^n \rho_\nu(x)= (-1)^n \rho_{\nu-n} (x),\quad n \in \mathbb{N}_0.\eqno(1.18)$$
In the meantime, the Mellin-Barnes integral (1.5) and reduction formula for the gamma-function yield 

$$\rho_{\nu+1} (x) = \frac{1}{2\pi i} \int_{\gamma-i\infty}^{\gamma+i\infty} \Gamma(\nu+s+1) \Gamma (s) x^{-s} ds$$

$$=  \frac{1}{2\pi i} \int_{\gamma-i\infty}^{\gamma+i\infty} \Gamma(\nu+s) (\nu+s) \Gamma (s) x^{-s} ds=  \nu \rho_\nu(x)$$

$$+  \frac{1}{2\pi i} \int_{\gamma-i\infty}^{\gamma+i\infty} \Gamma(\nu+s) \Gamma (s+1) x^{-s} ds$$

$$=   \nu \rho_\nu(x)+ x \rho_{\nu-1} (x).$$
Hence we deduce the following recurrence relation for the scaled Macdonald functions

$$\rho_{\nu+1} (x) =    \nu \rho_\nu(x)+ x \rho_{\nu-1} (x),\quad \nu \in \mathbb{R}.\eqno(1.19)$$
In the operator form it can be written as follows

$$\rho_{\nu+1} (x) = \left( \nu - xD \right) \rho_\nu(x).\eqno(1.20)$$
Further, recalling the definition of the operator $\theta$, identities (1.10)  and Rodrigues formula for the associated Laguerre polynomials,  we obtain

$$\theta^n \{ x^\nu e^{-x} \} =  n! x^{n+\nu} e^{-x} L_n^\nu(x),\  n \in \mathbb{N}_0.\eqno(1.21)$$
{\bf Lemma 1.}  {\it An arbitrary polynomial $f_n(x) = \sum_{k=0}^n f_{n,k} x^k$ of degree at most $n$ has the following integral representation

$$ f_n(x)= {1\over \rho_\nu(x) }  \int_0^\infty t^{\nu-1} e^{-t -x/t  }  q^\nu_{2n}(t) dt,\quad x >0,\eqno(1.22)$$ 
where $q^\nu_{2n}(x)$ of degree at most $2n$ is the associated polynomial given by the formula }

$$  q^\nu_{2n}(x)= \sum_{k=0}^n f_{n,k} (-1)^k k! x^k L_k^\nu(x).\eqno(1.23)$$

\begin{proof} Indeed, using (1.21),  (1.22),  we write the operator equality

$$f_n(-\theta) \left\{ x^\nu e^{-x}\right\} =  x^\nu e^{-x}  \sum_{k=0}^n f_{n,k} (-1)^k k! x^k L_k^\nu(x) = x^\nu e^{-x} q^\nu_{2n}(x).\eqno(1.24)$$
Then we plug in the left-hand side of the first equality in (1.24) into the right-hand side of (1.22) and integrate by parts,  eliminating the integrated terms, to find

$$\int_0^\infty t^{-1} e^{-x/t}  f_n(-\theta) \left\{ t^\nu e^{-t} \right\} dt  = \int_0^\infty f_n(\theta) \left\{ t^{-1} e^{-x/t} \right\} 
 t^\nu e^{-t} dt .$$
But

$$\theta^k \left\{ t^{-1} e^{-x/t} \right\}  = \left( t D t\right)^k \left\{ t^{-1} e^{-x/t} \right\} = x^k  t^{-1}  e^{-x/t}.$$
Hence, appealing to (1.6), we establish (1.22) and complete the proof of Lemma 1. 
\end{proof}

Besides, the representation (1.6) implies an immediate

{\bf Corollary 1}. { \it Assuming the explicit expression of the polynomial   $q^\nu_{2n}(x)= \sum_{k=0}^{2n} q_{n,k} x^k$, formula $(1.22)$ reads}

$$  f_n(x)= {1\over \rho_\nu(x) }    \sum_{k=0}^{2n} q_{n,k}\  \rho_{\nu+k}(x),\quad x >0.\eqno(1.25)$$

{ \bf Lemma 2} \ \cite{YAP}. {\it Let $n, m \in \mathbb{N},\nu \ge 0,\   f_n(x),\ g_{m}(x)$ be polynomials of  degree at most $n,\ m$, respectively.  Let 
 
 $$f_n(x) \rho_\nu(x)+ g_{m} (x) \rho_{\nu+1} (x) = 0\eqno(1.26)$$
 for all $x >0$. Then $f_n \equiv 0, \ g_{m} \equiv 0.$ }

{\bf Corollary 2}. { \it   Lemma 2 holds for the weight functions  $ \omega^\pm_\nu(x),\  \omega^\pm_{\nu+1} (x)$}.

{\bf Lemma 3}. {\it The weight functions $ \omega^\pm_\nu(x)$ satisfy the following second order differential equations, respectively, } 

$$x  [\omega^+_\nu(x) ]^{\prime\prime}+  (2x -\nu+1) [\omega^+_\nu(x)]^\prime + (x -\nu) \omega^+_\nu(x)=0,\eqno(1.27)$$

$$x^4   [\omega^-_\nu(x) ]^{\prime\prime} +  x^2 \left( (3-\nu) x  - 2   \right)  [\omega^-_\nu(x)]^\prime - (x-1)  \left(  x^2+\nu x +1\right)  \omega^-_\nu(x) =0.\eqno(1.28)$$

\begin{proof} In fact, recalling (1.19) and using (1.18), we have

$$[x  \omega^+_\nu(x) ]^\prime =  {d\over dx} \left[ e^{-x} (\rho_{\nu+2} (x)- (\nu+1)\rho_{\nu+1}(x)) \right]$$

$$= - (x-\nu-1)  \omega^+_\nu(x) -  e^{-x} \rho_{\nu+1} (x).$$
Hence with one more differentiation it gives

$$[x  \omega^+_\nu(x) ]^{\prime\prime} =  - (x -\nu-1) [\omega^+_\nu(x)]^\prime  + e^{-x} \rho_{\nu+1} (x).$$
Therefore summing these equations and fulfilling the differentiation, we arrive at (1.27).  Analogously, for the weight function  $ \omega^-_\nu(x)$ we find

$$[ x^2 \omega^-_\nu(x)]^\prime =   {d\over dx} \left[ e^{-1/x}  (\rho_{\nu+2} (x)- (\nu+1)\rho_{\nu+1}(x)) \right]
 $$

$$= \left(  (\nu+1) x + 1 \right) \omega^-_\nu(x) -  e^{-1/x}  \rho_{\nu+1} (x).$$
Hence

$$x^2 [ x^2 \omega^-_\nu(x)]^{\prime\prime} =    x^2 \left(  (\nu+1) x +  1  \right) [\omega^-_\nu(x)]^\prime +  x^2(\nu+1+x)  \omega^-_\nu(x)  -  e^{-1/x}  \rho_{\nu+1} (x),$$
and after subtraction of one equation from another   we get

$$  x^2 [ x^2 \omega^-_\nu(x)]^{\prime\prime}  - [ x^2 \omega^-_\nu(x)]^\prime = x^2 \left(  (\nu+1) x +  1  \right) [\omega^-_\nu(x)]^\prime + (x-1) \left(  (\nu+2) x + x^2 +1\right) \omega^-_\nu(x),$$
which yields (1.28) after simplification.

\end{proof}

Further, recalling  equality (1.14) and differentiating $n$ times its both sides, we find

$$S_n^\nu(x)= {d^n\over dx^n} \left[ x^n  \rho_\nu(x) \right]  =    n! \int_0^\infty t^{\nu -1} e^{-t - x/t}  L_n^\nu(t) dt,\eqno(1.29)$$
where the differentiation under the integral sign is allowed via the absolute and uniform convergence by $x \ge x_0 >0$.   As it is known \cite{YAP},  the sequence $S_n^\nu(x)$ generates the multiple orthogonal polynomials associated with the weight functions $(\rho_\nu(x),\ \rho_{\nu+1}(x))$, giving the corresponding Rodrigues formulas and explicit recurrence relation.   It is natural to characterize  multiple orthogonal polynomials for the weights 
$ ( \omega^\pm_\nu(x),\   \omega^\pm_{\nu+1}(x))$. However, the same scheme of investigation cannot be applied due to the features of these weights, in particular, their differential properties.   We leave this open question for the forthcoming research.

\section{Orthogonal polynomials with the weight $\omega^+_\nu(x)$} 

In this section we will investigate the sequence $P_n^\mu(x)$ of orthogonal polynomials with respect to the weight function $w^+_\nu(x)= e^{-x} \rho_\nu(x),\ \nu > -1$. Namely, our goal will be to establish its properties, the explicit representation, 3-term recurrence relation and Rodrigues-type formula. To do this,  we employ (1.3), (1.6),  Viskov-type identities   (1.10) and equality (1.21).  Indeed, we have

$$ \int_{0}^{\infty}  P^\nu_{n} (x)  e^{-x} \rho_\nu(x) \  x^m dx =  (-1)^m m!  \int_0^\infty t^{\nu+m -1} e^{-t }  L_m^\nu(t) \int_{0}^{\infty}  P^\nu_{n} (x)  e^{-x(1+1/t)} dx dt$$

$$=  (-1)^m  \int_0^\infty \theta^m \left\{ t^{\nu} e^{-t }\right\} {1\over t}  \int_{0}^{\infty}  P^\nu_{n} (x)  e^{-x(1+1/t)} dx dt$$

$$=  (-1)^m  \int_0^\infty \theta^m \left\{ t^{\nu} e^{-t }\right\}  P^\nu_n \left(\theta\right)  \left\{ {1\over t}  \int_{0}^{\infty}  e^{-x(1+1/t)} dx\right \}  dt$$

$$=  (-1)^m  \int_0^\infty \theta^m \left\{ t^{\nu} e^{-t }\right\}  P^\nu_n \left(\theta\right)  \left\{ {1\over1+ t} \right \}  dt$$

$$=    \int_0^\infty  t^{\nu} e^{-t }  P^\nu_n \left(\theta\right) \theta^m  \left\{ {1\over1+ t} \right \}  dt,\eqno(2.1)$$
where the interchange of the order of integration and differentiation under the integral sign are  allowed by  Fubini's  theorem and the absolute and uniform convergence.   Therefore, as it follows from (1.3) the  sequence $P_n^\nu$  is compositionally orthogonal in the sense of Laguerre (cf. \cite{YAP}) relatively to the function $(t+1)^{-1}$, i.e.

$$  \int_0^\infty  t^{\nu} e^{-t }  P^\nu_n \left(\theta\right) \theta^m  \left\{ {1\over1+ t} \right \}  dt =0,\quad  m=0,1,\dots, n-1.\eqno(2.2)$$
But it is easily seen via the Viskov-type identity that 

$$ \theta^m  \left\{ {1\over1+ t} \right \}  = t^m D^m  \left\{ {t^m \over1+ t} \right \} = (-1)^m t^m D^m  \left\{ {1 \over1+ t} \right \} =  {t^m \over (1+ t)^{m+1} }.$$
 Consequently, owing to integration by parts and (1.24), orthogonality conditions  (2.2) are  equivalent to the following equalities in terms of coefficients $a_{n,k}$ of the sequence $P^\nu_n (x)= \sum_{k=0}^n a_{n,k} x^k$

$$\int_0^\infty    t^{\nu+m} e^{-t } \sum_{k=0}^n a_{n,k} (-1)^k k! \ t^k L_k^\nu(t)  {dt\over (1+t)^{m+1}} =0,\quad   m=0,1,\dots, n-1.\eqno(2.3)$$
Moreover,  appealing to the method of mathematical induction, the latter equalities (2.3) can be reduced to the following relations

$$\int_0^\infty    t^{\nu} e^{-t } \sum_{k=0}^n a_{n,k} (-1)^k k! \ t^k L_k^\nu(t)  {dt\over (1+t)^{m+1}} =0,\quad   m=0,1,\dots, n-1.\eqno(2.4)$$
Writing,

 $${1\over (1+t)^{m+1} }  = {1\over m!} \int_0^\infty e^{-(1+t)y} y^m dy,$$
we substitute it in the left-hand side of (2.4)  and change the order of integration by Fubini's theorem.  Then, calculating the inner integral with respect to $t$  via \cite{PrudnikovMarichev}, Vol. II, Entry 2.19.3.2 

$$  \int_0^\infty    t^{\nu+k} e^{-t(1+y) } L_k^\nu(t)  dt = {\Gamma(1+\nu+k) \over (1+y)^{1+\nu+k}} \ p_k^{(\nu,0)} \left({y-1\over y+1} \right),\eqno(2.5)$$
where $p_k^{(\nu,0)}(x),\ k \in \mathbb{N}_0$ are the Jacobi polynomials \cite{Chi}, we obtain
$$  \int_0^\infty  e^{-y} y^m   \sum_{k=0}^n {a_{n,k} \ (-1)^k\  k! \  \Gamma(1+\nu+k)\over   (1+y)^{1+\nu+k}} \ p_k^{(\nu,0)} \left({y-1\over y+1}\right)   dy =0,\quad    m=0,1,\dots, n-1.\eqno(2.6)$$
Meanwhile, using the definition of the Jacobi polynomials in terms of the Gauss hypergeometric function \cite{PrudnikovMarichev}, Vol. II

$$p_k^{(\nu,0)} \left({y-1\over y+1} \right) = {(\nu+1)_k\over k!} {}_2F_1\left(-k,\ 1+\nu+k; \ 1+\nu;\ {1\over y+1} \right),\eqno(2.7)$$
equality (2.6) becomes
$$  \int_0^\infty  e^{-y} y^m   \sum_{k=0}^n  a_{n,k} \ (-1)^k  \Gamma ^2 (1+\nu+k) $$

$$\times  {1 \over   (1+y)^{1+\nu+k}} \  {}_2F_1\left(-k,\ 1+\nu+k; \ 1+\nu;\ {1\over y+1} \right)  dy =0,\quad    m=0,1,\dots, n-1.\eqno(2.8)$$
In the meantime, the sequence of functions 

$$F_n^\nu(y) = \sum_{k=0}^n  a_{n,k} \    {(-1)^k\   \Gamma ^2 (1+\nu+k) \over   (1+y)^{1+\nu+k}} \  {}_2F_1\left(-k,\ 1+\nu+k; \ 1+\nu;\ {1\over y+1} \right)\eqno(2.9)$$ 
can be expanded in terms of a series of the Laguerre polynomials  

$$F_n^\nu(y) = \sum_{k=0}^\infty c_{n,k} L_k(y),\eqno(2.10)$$
where

$$c_{n,k} =  \int_0^\infty  e^{-y}  L_k(y)  F_n^\nu(y) dy. \eqno(2.11)$$
Then, substituting the right-hand side of (2.10) into (2.8) and minding the orthogonality of the Laguerre polynomials,  we find

$$   \int_0^\infty  e^{-y} y^m   \sum_{k=0}^m  c_{n,k}\  L_k(y) dy =  0,\quad    m=0,1,\dots, n-1.\eqno(2.12)$$
But letting consecutively $m=0,1, \dots, n-1$ in (2.12), we obtain 

$$c_{n,m} =0, \quad m=0,1,\dots, n-1.\eqno(2.13)$$    
Returning to (2.9) and employing the definition of the Gauss hypergeometric function, we write

$$ F_n^\nu(y) = \sum_{k=0}^n  a_{n,k} \   (-1)^k \   \Gamma  (1+\nu+k)  \sum_{r=0}^k {(-k)_r\  \Gamma (1+\nu+k+r)\over 
(1+\nu)_r \ r! \ (1+y)^{1+\nu+k+r}}, \eqno(2.14)$$ 
where $(a)_k$ is the Pochhammer symbol \cite{Bateman}, Vol. I. Hence,  substituting the right-hand side of (2.14) in (2.11) and changing the order of integration and summation, we get 

$$c_{n,m} =  \sum_{k=0}^n     a_{n,k} \   (-1)^k \   \Gamma  (1+\nu+k)  \sum_{r=0}^k {(-k)_r\  \Gamma (1+\nu+k+r)\over 
(1+\nu)_r \ r! } \int_0^\infty    { e^{-y}  L_m(y) \over (1+y)^{1+\nu+k+r}} dy.\eqno(2.15) $$
The latter integral is calculated in \cite{PrudnikovMarichev}, Vol. II, Entry 2.19.3.9 in terms of the confluent hypergeometric function \cite{Bateman}, Vol. I, namely,

$$  \int_0^\infty    { e^{-y}  L_m(y) \over (1+y)^{1+\nu+k+r}} dy =  {1 \over \nu+k+r} \ {}_1F_1 \left( 1+m;\ 1-\nu-k-r;\  1\right) $$

$$+  {(1+ \nu+k+r)_m \over m!} \  \Gamma(-\nu-k-r) \ {}_1F_1 \left( 1+\nu+k+r+m;\ 1+\nu+k+r;\  1\right).\eqno(2.16)$$
Hence, substituting the value of the integral (2.16) in (2.15), we find

$$ c_{n,m} =  \sum_{k=0}^n     a_{n,k}\  d_{k,m},\eqno(2.17)$$
where, employing the reduction  and reflection formulas for the gamma-function \cite{Bateman}, Vol. I,

$$d_{k,m} =       \sum_{r=0}^k {(-k)_r \  (1+\nu)_k \over 
(1+\nu)_r \  r! }  \left[ (-1)^k \Gamma (\nu+k+r) \Gamma(1+\nu)  \ {}_1F_1 \left( 1+m;\ 1-\nu-k-r;\  1\right) \right.$$

$$\left. +   (-1)^{r}  {(1+ \nu+k+r)_m  \over  \ m!} \   \Gamma(-\nu) \Gamma^2(1+\nu)  {}_1F_1 \left( 1+\nu+k+r+m;\ 1+\nu+k+r;\  1\right) \right].\eqno(2.18)$$
Meanwhile,  the relation (7.11.1.21) in \cite{PrudnikovMarichev}, Vol. III says  

$$ \Gamma(\nu+k+r) \  {}_1F_1 \left( 1+m;\ 1-\nu-k-r;\  1\right) =   \Gamma (1+\nu+k+r+m )  \  \Psi \left( 1+m, \  1-\nu-k-r; 1 \right) $$

$$-  \Gamma (1+\nu) \Gamma(-\nu) (-1)^{k+r}  {(1+\nu+k+r)_m  \over m!} \  {}_1F_1 \left(  1+\nu+k+r+m;\ 1+\nu+k+r;\  1\right).\eqno(2.19)$$
Therefore,

$$ d_{k,m} =      \Gamma(1+\nu)    \sum_{r=0}^k  (-1)^{k+r}   \binom{k} {r} \  {  (1+\nu)_k \over (1+\nu)_r } \  \Gamma (1+\nu+k+r+m)   \Psi \left( 1+m, \  1-\nu-k-r; 1 \right) .\eqno(2.20)$$
In particular, we have

$$d_{0,m} =    \Gamma(1+\nu) \Gamma (1+\nu+m)   \Psi \left( 1+m, \  1-\nu; 1 \right).\eqno(2.21)$$
Hence, returning to (2.13), we observe that explicit values of the coefficients $a_{n,k},\ k =1,2,\dots, n$ can be expressed via Cramer's rule in terms of the free coefficient $a_{n,0}$ as follows  

$$  a_{n,k} = -    a_{n,0}\  { D_{n,k}  \over D_{n}},\quad   k = 1,\dots, n,\eqno(2.22)$$
where

$$ D_{n}= \begin{vmatrix}  

d_{1, 0}  & d_{2, 0}&  \dots&  \dots&   d_{n, 0} \\

d_{1, 1} &  \dots &   \dots&    \dots&    d_{n, 1} \\
 
\dots&  \dots &   \dots&   \dots&   \dots \\
 
 \vdots  &   \ddots  &  \ddots & \ddots&   \vdots\\
  
 d_{1, n-1} &   \dots&  \dots&   \dots&    d_{n, n-1}\\

 \end{vmatrix},\eqno(2.23)$$

$$ D_{n,k}= \begin{vmatrix}  

d_{1, 0}  & \dots & d_{k-1, 0}&  d_{0,0} &   d_{k+1, 0} & \dots&   d_{n, 0}  \\

d_{1, 1} &  \dots & d_{k-1, 1}&  d_{0,1} &   d_{k+1, 1} &  \dots&   d_{n, 1}  \\

 \dots&  \dots &   \dots&   \dots&  \dots&   \dots&  \dots \\
 
 \vdots  &   \ddots  &  \ddots & \vdots&  \ddots  &  \ddots & \vdots\\
  
 \vdots  &   \ddots  &  \ddots & \vdots&  \ddots  &  \ddots & \vdots\\

 \vdots  &   \ddots  &  \ddots & \vdots&  \ddots  &  \ddots & \vdots\\

d_{1, n-1} &  \dots & d_{k-1, n-1}&  d_{0,n-1} &   d_{k+1, 1} &  \dots&   d_{n, n-1}  \\

 \end{vmatrix}.\eqno(2.24)$$
The free  coefficient can be determined, in turn,  from the orthogonality conditions (1.1), (1.3), which imply the formula

$$\int_{0}^{\infty}  P^\nu_{n} (x)  e^{-x} \rho_\nu(x) \  x^n dx = {1\over a_{n,n}}.\eqno(2.25)$$ 
Therefore from (2.21), (2.22), (1.8) and the property for the Tricomi function (see \cite{PrudnikovMarichev}, Vol. III, Entry 7.11.4.2 ) we derive

$$ {1\over a_{n,n}} =  -   { a_{n,0} \over D_n \Gamma (1+\nu) } \sum_{k=0}^n  D_{n,k} (n+k)! \  d_{0,n+k},\eqno(2.26)$$
where $D_{n,0}\equiv - D_n.$  Hence

$$  a_{n,0}  =    \pm \  { D_n  \over [ D_{n,n} \Gamma(1+\nu) ]^{1/2} }\left[ \sum_{k=0}^n  D_{n,k} (n+k)! \ d_{0,n+k} \right]^{-1/2},\eqno(2.27)$$
where the sign can be  chosen  accordingly, making positive expressions under the square roots.   Assuming also the positivity of the leading coefficient $a_{n,n}$ we have its value, correspondingly,

$$   a_{n,n}  =    \mp \   [ D_{n,n}]^{1/2} \left[ {1\over \Gamma(1+\nu)} \sum_{k=0}^n  D_{n,k}  (n+k)! \ d_{0,n+k} \right]^{-1/2}.\eqno(2.28)$$
Furthermore, the sequence $ P_n^{\nu} (x)$ satisfies the 3-term recurrence relation in the form

$$x  P_n^{\nu} (x) = A_{n+1}  P_{n+1}^{\nu} (x) + B_n  P_n^{\nu} (x) + A_n  P_{n-1}^{\nu} (x),\eqno(2.29)$$
where $P_{-1}^{\nu} (x)\equiv 0$ and

$$A_{n+1}= {a_n\over a_{n+1} }, \quad B_{n}= {b_n\over a_n} -  {b_{n+1}\over a_{n+1}},\quad  a_n\equiv a_{n,n},\  b_n\equiv a_{n,n-1}.\eqno(2.30)$$
We summarize these results by the following 

{\bf Theorem 1}.  {\it Let $\nu > -1$. The sequence of orthogonal polynomials $P_n^{\nu} (x)$  can be expressed explicitly, where the coefficients $a_{n,k},\ k=1,2,\dots, n$ are calculated by formulas $(2.22)$ and the free term $a_{n,0}$ is defined by the equality $(2.27)$. Besides,  it satisfies the  3-term recurrence relation $(2.29)$, where }

$$A_{n+1}= { a_{n,0} \  D_{n+1} \ D_{n,n}   \over a_{n+1,0}\  D_n \ D_{n+1,n+1}},  \quad  B_n= { D_{n,n-1}  \over D_{n,n}} -  { D_{n+1,n}  \over D_{n+1,n+1}}.\eqno(2.31)$$

An analog of the Rodrigues formula for the sequence $ P^\nu_{n} (x)$ can be established, recalling  integral representation (1.22) of an arbitrary polynomial in terms of the associated polynomial of degree $2n$.  Consequently, using the differentiation under the integral sign owing to the absolute and uniform convergence with respect to $x \ge x_0 >0$, it is not difficult to obtain the following equality for the sequence  $ P^\nu_{n} (x)$

$$  P^\nu_{n} (x) = {e^x\  (-1)^n \over \rho_\nu(x) } {d^n\over dx^n}  \left( e^{-x}  \int_0^\infty t^{\nu+n-1} e^{-t -x/t  }  {q^\nu_{2n}(t)\over (1+t)^n } dt\right).\eqno(2.32)$$ 
In the meantime, expressing $q^\nu_{2n}(x)$ in terms of the associated Laguerre polynomials 

$$  q^\nu_{2n}(x)= \sum_{k=0}^{2n} f_{2n,k} L_k^\nu(x),\eqno(2.33)$$
we have, accordingly, the values of the coefficients, appealing to  (1.23), (2.22) and relation (2.19.14.15) in \cite{PrudnikovMarichev}, Vol. II, namely,

$$ f_{2n,k} =  {k!\over \Gamma(1+\nu+k)} \int_0^\infty t^\nu e^{-t} L_k^\nu(t)  q^\nu_{2n}(t) dt $$

$$= -  {k! \ a_{n,0} \over  D_n\ \Gamma(1+\nu+k)} \sum_{r=0}^n D_{n,r} (-1)^r r!  \int_0^\infty t^{\nu+r}  e^{-t} L_k^\nu(t)  L_r^\nu(t) dt $$

$$= -  { a_{n,0} \over  D_n} \sum_{r=0}^n D_{n,r} \  r! \ (1+\nu)_r  \ {}_3F_2 \left(-k,\ 1+\nu+r,\ 1+r;\ 1+\nu,\ 1;\ 1 \right). $$
Therefore

$$ f_{2n,k} = -  { a_{n,0} \over  D_n} \sum_{r=0}^n D_{n,r} \  r! \ (1+\nu)_r  \ {}_3F_2 \left(-k,\ 1+\nu+r,\ 1+r;\ 1+\nu,\ 1;\ 1 \right),\eqno(2.34) $$
where the values of the hypergeometric function can be simplified via relations (7.4.4; 90,91,92,93 ) in \cite{PrudnikovMarichev}, Vol. III. Precisely, we get for $k =0,1\dots, n$

$$ {}_3F_2 \left(-k,\ 1+\nu+r,\ 1+r;\ 1+\nu,\ 1;\ 1 \right) = 0,\quad k > 2r,\eqno(2.35)$$

$$ (1+\nu)_r\  {}_3F_2 \left(-k,\ 1+\nu+r,\ 1+r;\ 1+\nu,\ 1;\ 1 \right) = {  (2r) ! \over  r!  },\quad k=2r,$$ 

$$ (1+\nu)_r\  {}_3F_2 \left(-k,\ 1+\nu+r,\ 1+r;\ 1+\nu,\ 1;\ 1 \right) = {  (2r) ! (2r(\nu+2r) - r(\nu+r) ) \over  r! } \quad k=2r-1, $$

$$ (1+\nu)_r\  {}_3F_2 \left(-k,\ 1+\nu+r,\ 1+r;\ 1+\nu,\ 1;\ 1 \right) = {  (2r) ! \over 2\  r! }
\left( 2 r^2 (2r+\nu-1)(2r-1) \right.$$

$$\left. +  r(r-1)(r+\nu-1)(r+\nu ) \right),\quad k=2(r-1).$$ 
Substituting (2.34) in (2.33) and the result in (2.32), we have

$$  P^\nu_{n} (x) =  {e^x\  (-1)^{n+1}    a_{n,0}  \over  D_n\ \rho_\nu(x) } {d^n\over dx^n}  \sum_{k=0}^{2n}  \sum_{r=0}^n D_{n,r} \  r! \ (1+\nu)_r  \ {}_3F_2 \left(-k,\ 1+\nu+r,\ 1+r;\ 1+\nu,\ 1;\ 1 \right) $$

$$\times  e^{-x}  \int_0^\infty t^{\nu+n-1} e^{-t -x/t  }  {L_k^\nu(t)\over (1+t)^n } dt.\eqno(2.36)$$ 
The integral in (2.36) can be treated as follows. We write

$${t^n\over (1+t)^n} = {1\over n!} \int_0^\infty e^{ - (1+1/t) u} u^{n-1} du.$$
Hence, substituting the latter expression in (2.36), we change the order of integration by Fubini's theorem and recall (1.14) to obtain 

$$ \int_0^\infty t^{\nu+n-1} e^{-t -x/t  }  {L_k^\nu(t)\over (1+t)^n } dt =   {1\over n!} \int_0^\infty e^{ - u} u^{n-1}  \int_0^\infty t^{\nu-1} e^{-t -(x+u)/t  }  L_k^\nu(t) dt du $$

$$=   {(-1)^k \over n!} \int_0^\infty e^{ - u} u^{n-1}  {d^k\over dx^k } \int_0^\infty t^{\nu+k -1} e^{-t -(x+u) /t  }  L_k^\nu(t) dt du $$

$$=   {1 \over n!\  k! }  {d^k\over dx^k } \int_0^\infty e^{ - u} u^{n-1} \  (x+u)^k \rho_{\nu}(x+u)  du $$

$$=   {1 \over n!\  k! }  {d^k\over dx^k } \  e^x \int_x^\infty  (y-x)^{n-1} \   y^k  e^{ - y} \rho_{\nu}(y)  dy, $$
where the differentiation under the integral sign is allowed by virtue of the absolute and uniform convergence with respect to $x \ge x_0 >0$. Consequently, employing the definition (1.15) of the right-hand side Riemann- Liouville fractional integral and taking into account (2.35), equality (2.36) becomes 

$$  P^\nu_{n} (x) =  { (-1)^{n+1}    a_{n,0}  \over  D_n\  n }   \sum_{r=0}^{n} D_{n,r} \  r! \ (1+\nu)_r   \sum_{k=0}^{2r} {1\over k!}   \ {}_3F_2 \left(-k,\ 1+\nu+r,\ 1+r;\ 1+\nu,\ 1;\ 1 \right) $$

$$\times   {e^x\  \over \rho_\nu(x) } { d^n\over dx^n} \  e^{-x}  {d^k\over dx^k } \  e^x  I^n_ - \left(  x^k  e^{ - x} \rho_{\nu}(x) \right),\ n \in \mathbb{N}.\eqno(2.37)$$ 
The inicial polynomial $P_0^\nu(x)$ can be immediately obtained from (1.8) and (2.25),  and  we derive

$$P_0^\nu(x) = \left[  \Gamma(\nu+1) \Psi \left(1+\nu,\ 1+\nu; 1\right) \right]^{-1/2},\ \nu > -1.\eqno(2.38)$$
Concerning  the whole sequence  $P_n^\nu(x),\  n \in \mathbb{N}$, it can be established from

{\bf Theorem 2}. {\it Let $\nu > -1,\ n \in \mathbb{N}$. The sequence of orthogonal polynomials $P_n^\nu(x)$ satisfies the  Rodrigues-type  formula $(2.37)$, where the value $a_{n,0}$ is defined by formula $(2.27)$ and $D_n,\ D_{n,r}$ by $(2.23), (2.24)$, respectively.}

Further, property (1.19) and the orthogonality conditions (1.3) allow us to calculate the following values 

$$V_{n,\nu} =  \int_0^\infty  \left[ P_n^{\nu} (x) \right]^2  e^{-x} \rho_{\nu+1} (x) dx.\eqno(2.39)$$
In fact, using (1.19) and (1.1), we have

$$ V_{n,\nu} =  \nu +   \int_0^\infty  \left[ P_n^{\nu} (x) \right]^2  e^{-x} \rho_{\nu-1} (x) \ x dx.$$
Then integrating by parts in the latter integral and employing (2.25), (2.29), we get 

$$  V_{n,\nu} =  \nu +   1- \int_0^\infty  \left[ P_n^{\nu} (x) \right]^2  e^{-x} \rho_{\nu} (x) \ x dx + 2 \int_0^\infty   P_n^{\nu} (x) {d\over dx} \left[   P_n^{\nu} (x) \right]  e^{-x} \rho_{\nu} (x) \ x dx$$

$$=   \nu +   1- B_n  + 2 n ,$$
i.e.

$$ V_{n,\nu} =  \nu + 2n + 1- B_n.\eqno(2.40)$$
Finally, in this section we establish  the generating function for  polynomials $P_n^{\nu} (x)$, which is  defined as usually  by the equality

$$G(x,z) = \sum_{n=0}^\infty P_n^{\nu} (x) {z^n\over n!} ,\quad  x >0,\  z \in \mathbb{C},\eqno(2.41)$$
where $|z| < h_x$ and $h_x >0$ is a  convergence radius of the power series. To do this,  we employ (1.22), (1.29) and (2.33), having  the following equality from (2.41) 

$$  G(x,z) =  {1 \over \rho_\nu(x) } \sum_{n=0}^\infty  {z^n\over n!}  \sum_{k=0}^{2n} {f_{2n,k} \over k! }   {d^{k}\over dx^{k}} \left[  x^{k}  \rho_\nu(x) \right]=  {1\over \rho_\nu(x) } \sum_{n=0}^\infty  {z^n\over n!}  \sum_{k=0}^{2n}  f_{2n,k}   \sum_{j=0}^k  {(-1)^j\over j! (k-j)!} \  x^j \rho_{\nu-j}(x) . $$
Meanwhile, the product  $x^j\rho_{\nu-j}(x)$ is expressed in \cite{Cous}  as follows 

$$x^{ j}  \rho_{\nu-j}(x) = x^{ j/2} r_j(2\sqrt x; \nu ) \rho_{\nu}(x) +  x^{ (j-1)/2} r_{j-1} (2\sqrt x; \nu -1) \rho_{\nu+1}(x), \quad  j \in {\mathbb N}_0,$$
where $r_{-1}(z;\nu)=0$, 

$$ x^{ j/2} r_j(2\sqrt x; \nu ) = (-1)^j \sum_{i=0}^{[j/2]}   (\nu+i-j+1)_{j-2i} (j-2i+1)_i  {x^i\over i!}.\eqno(2.42)$$
Therefore this leads to the final  expression of  the generating function for the Prudnikov-type sequence $P_n^{\nu} (x)$, namely,

$$  G(x,z) =   \sum_{n=0}^\infty  \sum_{k=0}^{2n} \sum_{j=0}^k { (-1)^j \ f_{2n,k} \over n!\ j! (k-j)! }   x^{j/2}  r_j(2\sqrt x; \nu ) z^n $$

$$+  {\rho_{\nu+1}(x) \over \rho_\nu(x) } \sum_{n=0}^\infty \sum_{k=0}^{2n} \sum_{j=0}^k  { (-1)^j f_{2n,k} \over n!\ j! (k-j)! }  x^{(j-1)/2} r_{j-1}(2\sqrt x; \nu-1 ) z^n,\eqno(2.43) $$
where coefficients $f_{2n,k}$ are defined by (2.34).

\section{Orthogonal polynomials with the weight $w_\nu^-(x)$}

In the last section we will characterize the sequence $Q^\nu_{n} (x)$ of orthogonal polynomials with respect to the weight function  $w^-_\nu(x)= e^{- 1/x}  x^{-1} \rho_\nu(x),\ \nu  >0$.  Indeed, recalling (1.4), (1.6), (1.14), we treat its left-hand side similar to (2.1), having the chain of equalities

$$ \int_{0}^{\infty}  Q^\nu_{n} (x)  e^{- 1/x} \rho_\nu(x) \  x^{m-1} dx =  (-1)^m m!  \int_0^\infty t^{\nu+m -1} e^{-t }  L_m^\nu(t) \int_{0}^{\infty}  Q^\nu_{n} (x)  \ e^{- 1/x - x/t}\  {dx dt \over x} $$

$$=  (-1)^m  \int_0^\infty \theta^m \left\{ t^{\nu} e^{-t }\right\} {1\over t}  \int_{0}^{\infty}  Q^\nu_{n} (x) \  e^{- 1/x- x/t } \ {dx dt \over x} $$

$$=  (-1)^m  \int_0^\infty \theta^m \left\{ t^{\nu} e^{-t }\right\}  Q^\nu_n \left(\theta\right)  \left\{ {1\over t}  \int_{0}^{\infty}  e^{- 1/x - x/t}\  {dx\over x}\right \}  dt$$

$$=  (-1)^m  \int_0^\infty \theta^m \left\{ t^{\nu} e^{-t }\right\}  Q^\nu_n \left(\theta\right)  \left\{  {1\over t} \ \rho_0\left({1\over t}\right) \right \}  dt$$

$$=    \int_0^\infty  t^{\nu} e^{-t }  Q^\nu_n \left(\theta\right) \theta^m  \left\{  {1\over t} \ \rho_0\left({1\over t}\right) \right \}  dt.\eqno(3.1)$$
Consequently, conditions (1.4) imply that the  sequence $Q_n^\nu$  is compositionally orthogonal in the sense of Laguerre  relatively to the function $  1/t  \ \rho_0\left( 1/t \right)$, i.e.

$$  \int_0^\infty  t^{\nu} e^{-t }  Q^\nu_n \left(\theta\right) \theta^m  \left\{   {1\over t} \  \rho_0\left({1\over t} \right)\right \}  dt =0,\quad  m=0,1,\dots, n-1.\eqno(3.2)$$
In the meantime, employing equality (1.18), we obtain  

$$ \theta^m  \left\{  {1\over t} \  \rho_0\left({1\over t} \right) \right \}  =  \theta^{m-1}  \left\{  {1\over t} \  \rho_{-1} \left({1\over t} \right) \right \} = \dots =
  {1\over t} \  \rho_{-m} \left({1\over t} \right).\eqno(3.3)$$
Therefore after  integration by parts we find that  orthogonality conditions  (3.2) are  equivalent to the following equalities in terms of  coefficients  of the polynomial sequence $Q^\nu_n(x)$ (we keep the same notation  $a_{n,k}$ as in Section 2) 

$$\int_0^\infty    t^{\nu-1} e^{-t } \sum_{k=0}^n a_{n,k} (-1)^k k! \ t^k L_k^\nu(t)   \rho_{-m} \left({1\over t} \right) \ dt =0,\ \nu > 0, \quad   m=0,1,\dots, n-1.\eqno(3.4)$$
The integral in (3.4) can be  calculated via corrected formula (2.19.13.6) in \cite{PrudnikovMarichev}, Vol. II.  Namely, when $\alpha - \mu,\ \alpha, \ \mu \notin \mathbb{Z}$,  we obtain

$$ k! \int_0^\infty    t^{\alpha -1} e^{-t } L_k^\nu(t)   \rho_{\mu} \left({1\over t} \right) \ dt=   (1-\alpha+\mu+\nu)_k \ \Gamma( \alpha-\mu)\  \Gamma (-\mu) $$

$$\times \ {}_1F_3 \left(1-\alpha+\mu+\nu +k; \ 1+\mu, \ 1-\alpha+\mu,\  1-\alpha+\mu+\nu;\ -1\right) $$

$$+  (1-\alpha+\nu)_k \ \Gamma( \alpha)\  \Gamma (\mu)  \ {}_1F_3 \left(1-\alpha+\nu +k; \ 1-\mu, \ 1-\alpha,\  1-\alpha+\nu;\ -1\right) $$

$$+  (1+\nu)_k \ \Gamma( -\alpha)\  \Gamma (-\alpha + \mu)  \ {}_1F_3 \left(1+\nu +k; \ 1+\alpha, \ 1+\alpha- \mu,\  1+\nu;\ -1\right).\eqno(3.5) $$
Therefore, when $\nu- \mu,\ \nu, \ \mu \notin \mathbb{Z}$, we have 

$$ k! \int_0^\infty    t^{\nu+k -1} e^{-t } L_k^\nu(t)   \rho_{\mu} \left({1\over t} \right) \ dt$$

$$=   {\Gamma( \nu+k-\mu)\over \Gamma(1+\mu-k)} \   \Gamma (1+\mu)\ \Gamma (-\mu)  \ {}_0F_2 \left( 1-\nu-k+\mu,\  1-k+\mu;\ -1\right) $$

$$+  { \Gamma( \nu) \  \Gamma (\mu) \Gamma(1-\mu)\over \Gamma(1-\mu+k) }   \ {}_0F_2 \left(1-\nu,  \ 1-\mu+k ;\ -1\right) $$

$$+  (1+\nu)_k  \Gamma( -\nu-k)\  \Gamma (-\nu-k + \mu)  \ {}_0F_2 \left( 1+\nu+k- \mu,\  1+\nu;\ -1\right).\eqno(3.6) $$
Hence for $\nu \notin \mathbb{N} $

$$f_{k,m} = (-1)^k k!  \int_0^\infty    t^{\nu+k-1} e^{-t }  L_k^\nu(t)   \rho_{-m} \left({1\over t} \right) \ dt $$

$$= (-1)^k k! \lim_{\mu\to -m}  \int_0^\infty    t^{\nu+k-1} e^{-t }  L_k^\nu(t)   \rho_{\mu} \left({1\over t} \right) \ dt $$

$$= \lim_{\varepsilon \to 0} \left[  { (-1)^{k+m+1} \Gamma( \nu+k +m-\varepsilon)\over \varepsilon \Gamma(1+\varepsilon -m -k)} \   {}_0F_2 \left( 1-\nu-k-m+\varepsilon,\  1-k-m+\varepsilon;\ -1\right) \right.$$

$$+   { (-1)^{m+k} \ \Gamma( \nu) \over \varepsilon \Gamma(1+k+m-\varepsilon) }  \ {}_0F_2 \left(1-\nu,  \ 1+k+m -\varepsilon ;\ -1\right) $$

$$\left.+ (-1)^k  (1+\nu)_k  \Gamma( -\nu-k)\  \Gamma (-\nu-k -m+ \varepsilon)  \ {}_0F_2 \left( 1+\nu+k+m-\varepsilon,\  1+\nu;\ -1\right)\right] $$

$$= \lim_{\varepsilon \to 0} \left[  -  {1\over \varepsilon} \  \Gamma( \nu -\varepsilon)\Gamma(1-\nu+\varepsilon) \  \left[  \sum_{r=0}^{k+m-1}  {(-1)^r\over  r!\ \Gamma(  1-\nu-k-m+r+\varepsilon )\ \Gamma( 1-k-m+r+ \varepsilon) } \right.\right.$$

$$\left. + \sum_{r=k+m}^\infty {(-1)^r\over  r!\ \Gamma(  1-\nu-k-m+r+\varepsilon )\ \Gamma( 1-k-m+r+ \varepsilon) } \right]$$

$$\left. +   { (-1)^{m+k}\   \Gamma( \nu) \over \varepsilon\  \Gamma(1+k+m-\varepsilon)}   \ {}_0F_2 \left(1-\nu,  \ 1+k+m -\varepsilon ;\ -1\right) \right]$$

$$+   \Gamma( -\nu)\  \Gamma (-\nu-k -m)  \ {}_0F_2 \left( 1+\nu+k+m,\  1+\nu;\ -1\right) $$

$$=  \lim_{\varepsilon \to 0} \left[   -  {1\over \varepsilon} \  \Gamma( \nu -\varepsilon)\Gamma(1-\nu+\varepsilon) \ \sum_{r=0}^{k+m-1}  {(-1)^r\over  r!\ \Gamma(  1-\nu-k-m+r+\varepsilon )\ \Gamma( 1-k-m+r+ \varepsilon) } \right.$$

$$ +  {(-1)^{k+m+1}\over \varepsilon} \  {\Gamma( \nu -\varepsilon) \Gamma( 1-\nu +\varepsilon)\over \Gamma(1-\nu)  (k+m)!} \  \sum_{r=0}^\infty {(-1)^r\over  r!  (1-\nu)_r (1+k+m)_r }$$

$$ +   { (-1)^{m+k}\   \Gamma( \nu) \over \varepsilon\  \Gamma(1+k+m-\varepsilon)}   \ {}_0F_2 \left(1-\nu,  \ 1+k+m -\varepsilon ;\ -1\right) $$

$$\left.  +   (-1)^{k+m} {\Gamma( \nu -\varepsilon) \Gamma( 1-\nu +\varepsilon) \over  \Gamma(1-\nu) (k+m)!}  \  \sum_{r=0}^\infty {(-1)^r  \left[ \psi( 1-\nu+r )+ \psi( 1+r) \right] \over  r!  (  1-\nu )_r (1+k+m)_r  }  \right]$$

$$+   \Gamma( -\nu)\  \Gamma (-\nu-k -m)  \ {}_0F_2 \left( 1+\nu+k+m,\  1+\nu;\ -1\right) $$

$$=   (-1)^{k+m+1}  \sum_{r=0}^{k+m-1}  {(-1)^r\   \Gamma(\nu+ r+1)\  r! \over (k+m-r-1)! }   $$

$$  +   (-1)^{k+m} {\Gamma( \nu )  \over   (k+m)!}  \  \sum_{r=0}^\infty {(-1)^r  \left[ \psi( 1-\nu+r )+ \psi( 1+r) \right] \over  r!  (  1-\nu )_r (1+k+m)_r  } $$

$$  +   (-1)^{k+m} {\Gamma( \nu )  \over   (k+m)!}  \  \sum_{r=0}^\infty {(-1)^r  \  \psi( 1+k+m+r)  \over  r!  (  1-\nu )_r (1+k+m)_r  } $$

$$ +  (-1)^{k+m+1} { \Gamma( \nu) \over  (k+m)!}  \left[  \psi( 1-\nu ) -\psi( \nu ) \right] \sum_{r=0}^\infty {(-1)^r\over  r!  (1-\nu)_r (1+k+m)_r }$$

$$+   \Gamma( -\nu)\  \Gamma (-\nu-k -m)  \ {}_0F_2 \left( 1+\nu+k+m,\  1+\nu;\ -1\right), $$
where $\psi(z)$ is the digamma-function or Euler's psi-function \cite{Bateman}, Vol. I. Thus we established the following formula for the coefficients 
$f_{k,m}$

$$ f_{k,m} =   (-1)^{k+m+1}  \sum_{r=0}^{k+m-1}  { (-1)^r\ \Gamma(\nu+ r+1)\  r! \over (k+m-r-1)! }   $$

$$+  (-1)^{k+m} {\Gamma( \nu )\over (k+m)!}  \  \sum_{r=0}^\infty {(-1)^r  \over  r!  (  1-\nu)_r (1+k+m)_r } $$

$$\times  \left[ \psi( 1-\nu+r )+ \psi( 1+r) +  \psi( 1+k+m+r) -  \psi( 1-\nu ) + \psi( \nu )\right] $$

$$+   \Gamma( -\nu)\  \Gamma (-\nu-k -m)  \ {}_0F_2 \left( 1+\nu+k+m,\  1+\nu;\ -1\right) , \quad  \nu \notin \mathbb{N} .\eqno(3.7)$$
In order to obtain the values of coefficients $f_{k,m}$ for positive integers $\nu = l$, we  write

$$f_{k,m} =  \lim_{\varepsilon \to 0} \left[  (-1)^{k+m+1}  \sum_{r=0}^{k+m-1}  { (-1)^r\ \Gamma(l+ r+1+\varepsilon)\  r! \over (k+m-r-1)! }  \right. $$

$$+  (-1)^{k+m} {\Gamma( l+\varepsilon )\over (k+m)!}  \  \sum_{r=0}^\infty {(-1)^r  \over  r!  (  1-l-\varepsilon)_r (1+k+m)_r } $$

$$\times  \left[ \psi( 1+r- l - \varepsilon )+ \psi( 1+r) +  \psi( 1+k+m+r) -  \psi( 1- l - \varepsilon ) + \psi( l+ \varepsilon )\right] $$

$$\left. +   \Gamma( - l - \varepsilon)\  \Gamma (-l -k -m-\varepsilon)  \ {}_0F_2 \left( 1+l +k+m+\varepsilon,\  1+l + \varepsilon;\ -1\right) \right]. $$
Hence, employing the familiar identity for the digamma-function \cite{Bateman}, Vol. I

$$\psi(z)- \psi(1-z) = - \pi \cot (\pi z),$$
it gives

$$f_{k,m} =  \lim_{\varepsilon \to 0} \left[  (-1)^{k+m+1}  \sum_{r=0}^{k+m-1}  { (-1)^r\ \Gamma(l+ r+1+\varepsilon)\  r! \over (k+m-r-1)! }  \right. $$

$$+  (-1)^{k+m} {\Gamma( l+\varepsilon )\over (k+m)!}  \  \sum_{r=0}^\infty {(-1)^r  \over  r!  (  1-l-\varepsilon)_r (1+k+m)_r } $$

$$\times  \left[ \psi( 1+r- l - \varepsilon )+ \psi( 1+r) +  \psi( 1+k+m+r) \right] $$

$$+ {(-1)^{k+m+l+1} \over \varepsilon^2 }  \left(1- {(\pi\varepsilon)^2\over 2} \right)  \sum_{r=0}^\infty {(-1)^r  \over  r! \  \Gamma( 1+r-l-\varepsilon) \Gamma (1+k+m+r) } $$

$$\left. +   {(-1)^{k+m} \over \varepsilon^2}  \sum_{r=0}^\infty {(-1)^r  \over  r! \  \Gamma( 1+l +k+m+ r+ \varepsilon) \Gamma (1+l + r+ \varepsilon) }  \right]$$
$$=    (-1)^{k+m+1}  \sum_{r=0}^{k+m-1}  { (-1)^r\  (l+r)! \  r! \over (k+m-r-1)! } + \lim_{\varepsilon \to 0} \left[  (-1)^{k+m}  \  \sum_{r=0}^{l-1}  { \Gamma (l-r+\varepsilon) \over  r! \  \Gamma (1+k+m+r) } \right.$$

$$\times  \left[ \psi(  l-r+ \varepsilon )+  {1\over \varepsilon}  + \psi( 1+r) +  \psi( 1+k+m+r) \right] $$

$$+ {(-1)^{k+m}\over \varepsilon}  \  \sum_{r=0}^{\infty}  {(-1)^r \left[ \psi( 1+r - \varepsilon )+ \psi( 1+r+l) +  \psi( 1+k+m+r+l) \right] \over  (r+l) ! \  \Gamma (  1+r-\varepsilon) \Gamma (1+k+m+r+l) } $$

$$+ {(-1)^{k+m+1} \over \varepsilon }   \sum_{r=0}^{l-1}  {   \Gamma( l-r+\varepsilon)  \over  r! \   \Gamma (1+k+m+r) }  +  {(-1)^{k+m+1} \over \varepsilon^2 }  \sum_{r=0}^{\infty}  {(-1)^{r}  [ 1+\varepsilon \psi(1+r) ]\over  (r+l)! \ r! \  \Gamma (1+k+m+r+l) } $$

$$\left. +   {(-1)^{k+m} \over \varepsilon^2}  \sum_{r=0}^\infty {(-1)^r [ 1- \varepsilon [ \psi( 1+l+r) + \psi( 1+k+m+r +l ) ] + \varepsilon^2 \psi( 1+l+r) \psi( 1+k+m+r +l )] \over  r! \  \Gamma( 1+l +k+m+ r) \Gamma (1+l + r) }  \right]$$

$$=    (-1)^{k+m+1}  \sum_{r=0}^{k+m-1}  { (-1)^r\  (l+r)! \  r! \over (k+m-r-1)! } +  (-1)^{k+m}  \  \sum_{r=0}^{l-1}  { (l-r-1)!  \over  r! \  \Gamma (1+k+m+r) } $$

$$\times  \left[ \psi(  l-r )  + \psi( 1+r) +  \psi( 1+k+m+r) \right] $$

$$ +   (-1)^{k+m}   \sum_{r=0}^\infty {(-1)^r \  \psi( 1+l+r) \psi( 1+k+m+r +l )] \over  r! \  \Gamma( 1+l +k+m+ r) \Gamma (1+l + r) }.$$
Consequently, we derive the equality 

$$ f_{k,m} =   (-1)^{k+m}  \left[  \sum_{r=0}^{l-1}  { (l-r-1)!  \over  r! \  (k+m+r) ! }    \left[ \psi(  l-r )  + \psi( 1+r) +  \psi( 1+k+m+r) \right] \right.$$

$$\left. - \sum_{r=0}^{k+m-1}  { (-1)^r\  (l+r)! \  r! \over (k+m-r-1)! }   +     \sum_{r=0}^\infty {(-1)^r \  \psi( 1+l+r) \psi( 1+k+m+r +l ) \over  r! \  ( l +k+m+ r) ! \  (l + r) ! } \right],\ l \in \mathbb{N}.\eqno(3.8)$$
Returning to (3.4),  we write it in the form

$$ \sum_{k=0}^n a_{n,k} f_{k,m}  =0,\quad   m=0,1,\dots, n-1,\eqno(3.9)$$
observing via formulas (3.7), (3.8) that coefficients $f_{k,m} \equiv d_{k+m}$, i.e.  these values depend on the sum of indices.  Therefore equalities (3.9) can be expressed as a linear system of $n$ equations with $n$ unknowns, if we assume for now that the leading coefficient $a_{n,n}\equiv a_n > 0$ is a given value. Moreover,  this linear system involves  the Hankel matrix $\{ d_{k+m} \},\ (k,m) \in (0,1,\dots, n-1) \times   (0,1,\dots, n-1) $. Precisely,  we have

$$\begin{pmatrix}

d_{0}  & d_{1}&   d_2 &  \dots&   d_{n-1} \\

d_{1}&  d_{2} &   d_3&    \dots&   d_{n} \\
 
d_{2}&  d_{3} &   \dots& \dots&  d_{n+1} \\
 
\vdots  &  \vdots&   \ddots& \ddots&    \vdots\\
 
 d_{n-1} &  d_{n} &   d_{n+1} &  \dots& d_{2(n-1)}\\

\end{pmatrix} 
\begin{pmatrix}

a_{n,0} \\
 
 a_{n,1} \\

   \vdots  \\
       
    \vdots  \\
        
    a_{n, n-1}

\end{pmatrix}=  - a_n \ 
\begin{pmatrix}

 d_n \\
 
  d_{n+1}  \\

 \vdots  \\
 
  \vdots  \\
  
  d_{2n-1}
\end{pmatrix}
 .\eqno(3.10)$$
Hence, appealing to the  Cramer rule, we find the values of coefficients $a_{n,k}$ in terms of the leading coefficient $a_{n}$ in the form  

$$  a_{n,k} = -    a_{n}\  { D_{n,k}  \over D_{n}},\quad   k = 0,\dots, n-1,\eqno(3.11)$$
where

$$ D_{n}= \begin{vmatrix}

d_{0}  & d_{1}&   d_2 &  \dots&   d_{n-1} \\

d_{1}&  d_{2} &   d_3&    \dots&   d_{n} \\
 
d_{2}&  d_{3} &   \dots& \dots&  d_{n+1} \\
 
\vdots  &  \vdots&   \ddots& \ddots&    \vdots\\
 
 d_{n-1} &  d_{n} &   d_{n+1} &  \dots& d_{2(n-1)}\\

\end{vmatrix},\eqno(3.12)$$

$$ D_{n,k}= \begin{vmatrix}  

d_{0}  & \dots& d_{k-1}&  d_{n} &   d_{k+1} & \dots&   d_{n-1}  \\

d_{1} &  \dots & d_{k}&  d_{n+1} &   d_{k+2} &  \dots&   d_{n}  \\

 d_2  &   \ddots & \vdots&  \vdots  &  \ddots & \vdots& d_{n+1}\\
  
 \vdots  &   \ddots  &  \ddots & \vdots&  \ddots  &  \ddots & \vdots\\

 \vdots  &   \ddots  &  \ddots & \vdots&  \ddots  &  \ddots & \vdots\\

d_{ n-1} &  \dots & d_{k+n-2}&  d_{2n-1} &   d_{k+n} &  \dots&   d_{2(n-1)}  \\

 \end{vmatrix}.\eqno(3.13)$$
Then the leading term $a_n$ can be defined, employing an analog of the equality (2.25) for the sequence  $Q^\nu_{n} (x)$

$$ \int_{0}^{\infty}  Q^\nu_{n} (x)  e^{- 1/x} \rho_\nu(x) \  x^{n-1} dx =  {1\over a_n}.$$
Indeed, we write from (3.11) and the latter equality

$$ {1\over a_n} = \sum_{k=0}^n a_{n,k}  \int_{0}^{\infty}    e^{- 1/x} \rho_\nu(x) \  x^{n+k-1} dx = - {a_n\over D_n} \sum_{k=0}^n D_{n,k}  \int_{0}^{\infty}    w^-_\nu(x) \  x^{n+k} dx,\eqno(3.14)$$
where $D_{n,n}= - D_n$ and  the $n+k$-th moment of the weight function $w^-_\nu(x)$ can be calculated via formula (1.9). In fact, it has for $\nu \notin \mathbb{N}$

$$\mu_{n+k}^\nu=  \int_0^\infty \omega^-_\nu(x) x^{n+k} dx =  \lim_{\varepsilon \to 0} \left[  {1\over \varepsilon} \ 
\sum_{r=0}^\infty  { (-1)^r \  \Gamma(\varepsilon+\nu) \Gamma(1-\varepsilon-\nu)  \over r! \ \Gamma(r+1-n-k-\varepsilon) \Gamma (r+1-n-k-\varepsilon- \nu) } \right.$$

$$\left. +  (-1)^{n+k}  \sum_{r=0}^\infty  { (-1)^r \  \Gamma(-\nu)\Gamma(1+\nu)  \Gamma(-\varepsilon-\nu)  \Gamma(1+\varepsilon +\nu)
 \over r! \ \Gamma(r+1+n+k+\varepsilon+\nu) \Gamma (r+1+ \nu) } \right.$$
 
 $$\left.  + { (-1)^{n+k+1} \over \varepsilon}  \sum_{r=0}^\infty  { (-1)^r \  \Gamma(\nu)\Gamma(1-\nu)   \over r! \ \Gamma(r+1+n+k+\varepsilon) \Gamma (r+1-\nu) } \right] $$

$$=    \sum_{r=0}^{n+k-1}   {  (-1)^r  \over r!}  \  (n+k-1-r)! \  \Gamma(n+k-r+\nu) $$

$$+  {(-1)^{n+k} \Gamma(\nu)  \over (n+k)!}   \sum_{r= 0}^ \infty  { (-1)^{r}\   [  \psi( r+1) +  \psi( r+1+n+k) + \psi( r+1-\nu) ]  \over r! \ (1+n+k)_r \  (1-\nu)_r } $$

$$+  {(-1)^{n+k}  \Gamma(\nu) \over (n+k)!} \left[ \psi(\nu)  -   \psi(1-\nu) \right]  \sum_{r= 0}^ \infty  { (-1)^{r}\    \over r! \ (1+n+k)_r \  (1-\nu)_r } $$

$$ +  (-1)^{n+k} [ \Gamma (-\nu) ]^2 \sum_{r=0}^\infty  { (-1)^r  \over r! \  (1+\nu)_{r+n+k}\  (1+ \nu)_r } .$$
 Therefore we get the formula
 
$$  \mu_{n+k}^\nu =   \sum_{r=0}^{n+k-1}   {  (-1)^r  \over r!}  \  (n+k-1-r)! \  \Gamma(n+k-r+\nu) $$

$$+  (-1)^{n+k} \Gamma(\nu)  \sum_{r= 0}^ \infty  { (-1)^{r}\   [  \psi( r+1) +  \psi( r+1+n+k) + \psi( r+1-\nu) -  \pi \cot(\pi\nu) ]  \over r! \ (r+n+k) ! \  (1-\nu)_r } $$

$$ +  {(-1)^{n+k} [\Gamma (-\nu)]^2  \over (1+\nu)_{n+k}}\   {}_0F_2 \left( 1+\nu +n+k,\  1+\nu;\ -1\right) ,\quad  \nu \notin \mathbb{N}.\eqno(3.15)$$
When $\nu \in \mathbb{N}$ it reads

$$ \mu_{n+k}^\nu =   \lim_{\varepsilon \to 0} \left[  \sum_{r=0}^{n+k-1}   {  (-1)^r  \over r!}  \  (n+k-1-r)! \  \Gamma(n+k-r+\nu+\varepsilon) \right.$$

$$ +  (-1)^{n+k+\nu} {\pi\over \sin(\pi\varepsilon)}   \sum_{r= 0}^ \infty  { (-1)^{r}\   [  \psi( r+1) +  \psi( r+1+n+k) + \psi( r+1-\nu-\varepsilon) - \pi\cot(\pi\varepsilon) ]  \over r! \ \Gamma (1+r+n+k) \  \Gamma(r+1-\nu-\varepsilon) } $$

$$\left. +  (-1)^{n+k} [\Gamma (-\nu-\varepsilon)]^2   \sum_{r=0}^\infty  { (-1)^r  \over r! \  (1+\nu+\varepsilon)_{r+n+k}\  (1+ \nu+\varepsilon)_r } \right]$$

$$=   \sum_{r=0}^{n+k-1}   {  (-1)^r  \over r!}  \  (n+k-1-r)! \  (n+k-r+\nu-1) ! $$

$$ +  \lim_{\varepsilon \to 0} \left[  {(-1)^{n+k}  \over (n+k)!}   \sum_{r= 0}^{\nu-1}  {  \psi( r+1) +  \psi( r+1+n+k) + \psi( r+1-\nu-\varepsilon) - \pi\cot(\pi\varepsilon)  \over r! \ (1+n+k)_r } \Gamma(\nu-r+\varepsilon) \right.$$

$$+ (-1)^{n+k} {\pi\over \sin(\pi\varepsilon)}   \sum_{r= 0}^ \infty  { (-1)^{r}\   [  \psi( r+\nu+1) +  \psi( r+1+n+k+\nu) + \psi( r+1-\varepsilon)- \pi\cot(\pi\varepsilon) ]  \over (r+\nu)! \ \Gamma (1+n+k+r+\nu) \  \Gamma(r+1-\varepsilon) } $$

$$\left. +  {(-1)^{n+k} \pi^2 \over \sin^2(\pi \varepsilon)}  \sum_{r=0}^\infty  { (-1)^r  \over r! \  \Gamma (1+\nu+r+n+k+\varepsilon)\  \Gamma (1+ \nu+r+\varepsilon) } \right]$$

$$=   \sum_{r=0}^{n+k-1}   {  (-1)^r  \over r!}  \  (n+k-1-r)! \  (n+k-r+\nu-1) ! $$

$$ +  {(-1)^{n+k}  \over (n+k)!}   \sum_{r= 0}^{\nu-1}  {  (\nu-r-1)!\ [ \psi( r+1) +  \psi( r+1+n+k) + \psi( \nu-r) ] \over r! \ (1+n+k)_r } $$

$$+  \lim_{\varepsilon \to 0} \left[  {(-1)^{n+k} \pi \over \sin (\pi\varepsilon)}     \sum_{r= 0}^ \infty  { (-1)^{r}\   [  \psi( r+\nu+1) +  \psi( r+1+n+k+\nu) + \psi( r+1-\varepsilon) - \pi\cot(\pi\varepsilon)  ]  \over (r+\nu)! \ \Gamma (1+n+k+r+\nu) \  \Gamma(r+1-\varepsilon) } \right.$$

$$\left. +  {(-1)^{n+k} \pi^2\over \sin^2(\pi\varepsilon)}  \sum_{r=0}^\infty  { (-1)^r  \over r! \  \Gamma (1+\nu+r+n+k+\varepsilon)\  \Gamma (1+ \nu+r+\varepsilon) } \right]$$
$$=   \sum_{r=0}^{n+k-1}   {  (-1)^r  \over r!}  \  (n+k-1-r)! \  (n+k-r+\nu-1) ! $$

$$ +  (-1)^{n+k}  \left[  \sum_{r= 0}^{\nu-1}  {  (\nu-r-1)!\ [ \psi( r+1) +  \psi( r+1+n+k) + \psi( \nu-r) ] \over r! \ (n+k+r)! } \right.$$

$$ +  \sum_{r= 0}^ \infty  { (-1)^{r}\  \psi(r+1)   [  \psi( r+\nu+1) +  \psi( r+1+n+k+\nu)  ]  \over  r !\ (r+\nu)! \  (n+k+r+\nu) !  } $$

$$+  \sum_{r=0}^\infty  { (-1)^r   \psi (1+\nu+r+n+k) \psi (1+\nu+r)  \over r! \  (\nu+r+n+k)! \ ( \nu+r)! } $$

$$+ {1\over 2}   \sum_{r=0}^\infty  { (-1)^r [ \psi^2 (r+1+\nu+n+k)  - \psi^\prime (r+1+\nu+n+k) ]   \over r! \  (\nu+r+n+k)! \ ( \nu+r)! } $$

$$\left. + {1\over 2}  \sum_{r=0}^\infty  { (-1)^r [ \psi^2 (r+1+\nu)  - \psi^\prime (r+1+\nu) ]   \over r! \  (\nu+r+n+k)! \ ( \nu+r)! } \right]$$

$$+   (-1)^{n+k}  \lim_{\varepsilon \to 0} \  \left[  {\pi \over \sin (\pi\varepsilon)}   \sum_{r= 0}^ \infty  { (-1)^{r}\  \psi( r+1-\varepsilon)  \over \Gamma(r+1-\varepsilon) (r+\nu)! \  (n+k+r+\nu)! } \right.$$

$$\left.  -  {\pi^2   \cot(\pi\varepsilon) \over  \sin (\pi\varepsilon)  }   \sum_{r= 0}^ \infty  { (-1)^{r}\   \over \Gamma(r+1-\varepsilon) (r+\nu)! \  (n+k+r+\nu)! } + {\pi^2\over \sin^2(\pi\varepsilon)}  \sum_{r=0}^\infty  { (-1)^r  \over r! \  (\nu+r+n+k)!  ( \nu+r)! }\right]$$

$$=   \sum_{r=0}^{n+k-1}   {  (-1)^r  \over r!}  \  (n+k-1-r)! \  (n+k-r+\nu-1) ! $$

$$ +  (-1)^{n+k}  \left[  \sum_{r= 0}^{\nu-1}  {  (\nu-r-1)!\ [ \psi( r+1) +  \psi( r+1+n+k) + \psi( \nu-r) ] \over r! \ (n+k+r)! } \right.$$

$$ + \sum_{r= 0}^ \infty  { (-1)^{r}\  \psi(r+1)   [  \psi( r+\nu+1) +  \psi( r+1+n+k+\nu)  ]  \over  r !\ (r+\nu)! \  (n+k+r+\nu) !  } $$

$$\left. +  \sum_{r=0}^\infty  { (-1)^r  [\psi (1+\nu+r+n+k) \psi (1+\nu+r) +\pi^2/2]  \over r! \  (\nu+r+n+k)! \ ( \nu+r)! } \right.$$

$$ + {1\over 2}  \sum_{r=0}^\infty  { (-1)^r [ \psi^2 (r+1)  - \psi^\prime (r+1) ]   \over r! \  (\nu+r+n+k)! \ ( \nu+r)! } $$

$$+ {1\over 2}   \sum_{r=0}^\infty  { (-1)^r [ \psi^2 (r+1+\nu+n+k)  - \psi^\prime (r+1+\nu+n+k) ]   \over r! \  (\nu+r+n+k)! \ ( \nu+r)! } $$

$$\left. + {1\over 2}  \sum_{r=0}^\infty  { (-1)^r [ \psi^2 (r+1+\nu)  - \psi^\prime (r+1+\nu) ]   \over r! \  (\nu+r+n+k)! \ ( \nu+r)! } \right]. $$
Thus after a slight simplification we have finally,

$$ \mu_{n+k}^\nu =  \sum_{r=0}^{n+k-1}   {  (-1)^r  \over r!}  \  (n+k-1-r)! \  (n+k-r+\nu-1) ! $$

$$ +  (-1)^{n+k}  \left[  \sum_{r= 0}^{\nu-1}  {  (\nu-r-1)!\ [ \psi( r+1) +  \psi( r+1+n+k) + \psi( \nu-r) ] \over r! \ (n+k+r)! } \right.$$

$$ + {1\over 2}  \sum_{r=0}^\infty  { (-1)^r [ \psi (r+1) + \psi (r+1+\nu) +  \psi (r+1+\nu+n+k)]^2    \over r! \  (\nu+r+n+k)! \ ( \nu+r)! } $$

$$\left.   -   {1\over 2}  \sum_{r=0}^\infty  { (-1)^r [  \psi^\prime (r+1) + \psi^\prime (r+1+\nu)+ \psi^\prime (r+1+\nu+n+k)  ]   \over r! \  (\nu+r+n+k)! \ ( \nu+r)! } \right]$$

$$ + {\pi^2\ (-1)^{n+k}  \over 2\  (\nu+n+k)! \  \nu ! }\  {}_0F_2 \left( 1+\nu+n+k,\ 1+\nu;\ -1\right),\  \nu \in \mathbb{N}.\eqno(3.16)$$
Taking into account the moment values (3.15), (3.16), we return to (3.14) to find the value of the leading term $a_n$ in the form

$$  a_n = \left[  - {1 \over D_n} \sum_{k=0}^n D_{n,k} \  \mu_{n+k}^\nu\right]^{- 1/2}. \eqno(3.17)$$
Hence we arrive at

{\bf Theorem 3}.  {\it Let $\nu > 0$. The sequence of orthogonal polynomials $Q_n^{\nu} (x)$  can be expressed explicitly, where the coefficients $a_{n,k},\ k=0,2,\dots, n-1$ are calculated by formulas $(3.11)$ and the leading  term $a_{n}$ is given by the equality $(3.17)$. Besides,  it satisfies the  3-term recurrence relation

$$x  Q_n^{\nu} (x) = A_{n+1}  Q_{n+1}^{\nu} (x) + B_n  Q_n^{\nu} (x) + A_n  Q_{n-1}^{\nu} (x),\quad n \in \mathbb{N}_0$$
where $Q_{-1}^{\nu} (x)\equiv 0$ and }

$$A_{n+1}= { a_{n}  \over a_{n+1}},  \quad  B_n= { D_{n+1,n}  \over D_{n+1} }-  { D_{n,n-1}  \over D_{n}}.$$

An analog of the Rodrigues formula for the sequence $ Q^\nu_{n} (x)$ can be obtained in a similar manner as in the previous section for the sequence $P_n^\nu$. In fact, recalling (1.15), (1.22), keeping the same notations  and taking into account the permission of the differentiation under the integral sign and the interchange of the order of integration, we write

$$ Q^\nu_n(x)= {e^{1/x} x \  (-1)^n \over \rho_\nu(x) (n-1)! }  {d^n\over dx^n} \int_0^\infty t^{\nu-1} e^{-t  }  q^\nu_{2n}(t) \int_x^\infty (y-x)^{n-1} e^{-1/y- y/t} {dy dt \over y} $$ 

$$= {e^{1/x} x \  (-1)^n \over \rho_\nu(x) (n-1)! }  {d^n\over dx^n} \int_0^\infty t^{\nu-1} e^{-t -x/t }  q^\nu_{2n}(t) \int_0^\infty y^{n-1} e^{-1/(y+x)- y/t} {dy dt \over y+x} .$$ 
But the inner integral with respect to $y$ can be calculated in terms of the series of Tricomi's functions. Namely, justifying the interchange of the order of summation and integration by the absolute and uniform convergence with respect to $x \ge x_0 >0$, we derive

$$ \int_0^\infty y^{n-1} e^{-1/(y+x)- y/t} {dy \over y+x}  = \sum_{k=0}^\infty {(-1)^k\over k!}  \int_0^\infty { y^{n-1} e^{- y/t} \over (y+x)^{k+1} }\ dy $$

$$= (n-1)!  \sum_{k=0}^\infty {(-1)^k x^{n-k-1} \over k!} \Psi \left( n,\  n-k;\  {x\over t} \right).$$
Therefore we obtain

$$ Q^\nu_n(x)= {e^{1/x} x \  (-1)^n \over \rho_\nu(x)  }  {d^n\over dx^n}   \sum_{k=0}^\infty {(-1)^k x^{n-k-1} \over k!} \int_0^\infty t^{\nu-1} e^{-t  }  q^\nu_{2n}(t)  \Psi \left( n,\  n-k;\  {x\over t} \right) dt.\eqno(3.18)$$
Meanwhile, the latter integral with respect to $t$ can be treated with the use of (1.23), (3.11), relations (8.4.33.3), (8.4.46.1) in \cite{PrudnikovMarichev}, Vol. III  and the Parseval equality for the Mellin transform \cite{Tit}.  In fact, we have

$$\int_0^\infty t^{\nu-1} e^{-t  }  q^\nu_{2n}(t)  \Psi \left( n,\  n-k;\  {x\over t} \right) dt =     { a_{n} \over D_{n}} \sum_{r=0}^n  (-1)^{r+1} r! \  D_{n,r} 
\int_0^\infty t^{\nu-1} e^{-t  }  L_r^\nu (t)  \Psi \left( n,\  n-k;\  {x\over t} \right) dt $$
and, accordingly, appealing to relation (2.16.3.14) in \cite{PrudnikovMarichev}, Vol. II

$$ \int_0^\infty t^{\nu-1} e^{-t  }  L_r^\nu (t)  \Psi \left( n,\  n-k;\  {x\over t} \right) dt$$

$$ = {1\over 2\pi i \ (n-1)!\ k! \ r! } \int_{\gamma -i \infty}^{\gamma +i\infty}
\frac{  \Gamma(s)\ \Gamma( s+\nu)   \Gamma(1+k-n+s)  \Gamma (1+r-s) \Gamma (n-s)}{ \Gamma(1-s) } x^{-s} ds $$

$$ = {1\over 2\pi i \ (n-1)!\ k! \ r! }   {d^r\over dx^r}   \int_{\gamma -i \infty}^{\gamma +i\infty}
 \Gamma(s) \Gamma( s+\nu)   \Gamma(1+k-n+s) \Gamma (n-s)\  x^{r-s} ds $$

$$ = {(-1)^{k+n+1} \over 2\pi i \ (n-1)!\ k! \ r! }   {d^r\over dx^r}  x^{r+k-n}  {d^k\over dx^k}   \int_{\gamma -i \infty}^{\gamma +i\infty}
 \Gamma^2 (s) \Gamma( s+\nu)   \Gamma(1-s) \  x^{n-1-s} ds $$

$$ = {(-1)^{k+n+1}  2^{2-\nu} \over \ (n-1)!\ k! \ r! }   {d^r\over dx^r}  x^{r+k-n}  {d^k\over dx^k} x^{n-1}  \int_0^\infty { t^{\nu+1} K_\nu (t) \over t^2 +4x} dt  $$

$$ = {(-1)^{k+n+1}  2^{2+\nu} \Gamma(1+\nu) \over \ (n-1)!\ k! \ r! }   {d^r\over dx^r}  x^{r+k-n}  {d^k\over dx^k} \left[ x^{n+\nu/2 -1}  S_{-\nu-1, \nu} (2\sqrt x) \right],$$
where $S_{\mu,\nu}(z)$ is the Lommel function \cite{Bateman}, Vol. II.  Hence combining with (3.18), we establish the following Rodrigues-type formula for the orthogonal sequence   $Q^\nu_n$

$$ Q^\nu_n(x)= {e^{1/x} x \  2^{2+\nu} \Gamma(1+\nu)  a_n \over \rho_\nu(x)\  D_n \ (n-1)!}    \sum_{k=0}^\infty \sum_{r=0}^n  {(-1)^{r}  \over  [ k!]^2 } \  D_{n,r} $$

$$\times  { d^n\over dx^n} \  x^{n-k-1} \  {d^r\over dx^r}  x^{r+k-n}  {d^k\over dx^k} \left[ x^{n+\nu/2 -1}  S_{-\nu-1, \nu} (2\sqrt x) \right].\eqno(3.19)$$

{\bf Theorem 4}.  {\it Let $\nu > 0$. The sequence of orthogonal polynomials $Q_n^{\nu} (x)$  can be generated via the Rodrigues-type formula $(3.19)$, where $a_n, D_n,\ D_{n,r} $ are calculated by formulas $(3.12), (3.13), (3.17)$}.

Now, recalling  (1.4),  (1.19), we calculate the values 

$$M_{n,\nu} =  \int_0^\infty  \left[ Q_n^{\nu} (x) \right]^2  e^{-1/x} \rho_{\nu+1} (x) dx.\eqno(3.20)$$
In fact, employing (1.19) and (1.2), we find

$$ M_{n,\nu} =  \nu \int_0^\infty  \left[ Q_n^{\nu} (x) \right]^2  e^{-1/x} \rho_{\nu} (x) dx +   \int_0^\infty  \left[ Q_n^{\nu} (x) \right]^2  e^{-1/x} \rho_{\nu-1} (x)  x dx.\eqno(3.21)$$
Then,   using the 3-term recurrence relation for the sequence $Q_n^\mu$,  we immediately obtain 

$$ \int_0^\infty  \left[ Q_n^{\nu} (x) \right]^2  e^{-1/x} \rho_{\nu} (x) dx = B_n.$$
Moreover, integrating by parts in the second integral in (3.21), it gives

$$\int_0^\infty  \left[ Q_n^{\nu} (x) \right]^2  e^{-1/x} \rho_{\nu-1} (x)  x dx = B_n - 1 +  2 \int_0^\infty   Q_n^{\nu} (x) {d\over dx} \left[   Q_n^{\nu} (x) \right]  e^{-1/x} \rho_{\nu} (x) x dx.$$
But employing again the 3-term recurrence relation for the sequence $Q_n^\mu$, it yields  $(b_n \equiv a_{n,n-1})$

$$\int_0^\infty   Q_n^{\nu} (x) {d\over dx} \left[   Q_n^{\nu} (x) \right]  e^{-1/x} \rho_{\nu} (x) x dx = (n-1) {b_n\over a_n} - n  {b_{n+1}\over a_{n+1} } =  
n  B_n -  {b_n\over a_n}.$$
Hence,  substituting these values in (3.21), we get, finally,  values (3.20) in the form

$$ M_{n,\nu} = ( \nu +1+2n) B_n  -1 - {2 b_n\over a_n}.\eqno(3.22)$$
The corresponding generating function for the sequence $Q_n^{\nu} (x)$

$$F(x,z) = \sum_{n=0}^\infty Q_n^{\nu} (x) {z^n\over n!} ,\quad  x >0,\  z \in \mathbb{C},$$
where $|z| < h_x$ can be obtained similarly to (2.43).  Indeed,  recalling (1.22), (1.29), (2.33),  (2.42),  we derive

$$  G(x,z) =   \sum_{n=0}^\infty  \sum_{k=0}^{2n} \sum_{j=0}^k { (-1)^j \ f_{2n,k} \over n!\ j! (k-j)! }   x^{j/2}  r_j(2\sqrt x; \nu ) z^n $$

$$+  {\rho_{\nu+1}(x) \over \rho_\nu(x) } \sum_{n=0}^\infty \sum_{k=0}^{2n} \sum_{j=0}^k  { (-1)^j f_{2n,k} \over n!\ j! (k-j)! }  x^{(j-1)/2} r_{j-1}(2\sqrt x; \nu-1 ) z^n,\eqno(3.23) $$
where for this case the coefficients $f_{2n,k}$ are defined by the formula (cf. (3.11), (3.12), (3.13), (3.17) )
$$ f_{2n,k} = -  { a_{n} \over  D_n} \sum_{r=0}^n D_{n,r} \  r! \ (1+\nu)_r  \ {}_3F_2 \left(-k,\ 1+\nu+r,\ 1+r;\ 1+\nu,\ 1;\ 1 \right).$$

\bibliographystyle{amsplain}

\end{document}